\documentclass{article}

\usepackage{amsmath,amssymb,amsthm}
\usepackage{graphicx}

\title{Optimal Control of Newton-Type\\
Problems of Minimal Resistance\footnote{This research was
partially presented at the Second Junior European Meeting
``Control Theory and Stabilization'', Torino, Italy, 3--5 December
2003. Research report CM04/I-01, Dep. Mathematics, Univ. Aveiro,
January 2004. Accepted (25-03-2004) for publication in the
\emph{Rendiconti del Seminario Matematico
dell'Universit\`{a} e del Politecnico di Torino}.}}

\author{Delfim F.~M.~Torres \\ \texttt{delfim@mat.ua.pt}
        \and
        Alexander Yu. Plakhov \\ \texttt{plakhov@mat.ua.pt}
        }

\date{Department of Mathematics\\
      University of Aveiro\\
      3810-193 Aveiro, Portugal}

\newtheorem{theorem}{Theorem}
\newtheorem{corollary}[theorem]{Corollary}
\newtheorem{proposition}[theorem]{Proposition}

\theoremstyle{remark}
\newtheorem{remark}[theorem]{Remark}

\theoremstyle{definition}
\newtheorem{definition}[theorem]{Definition}

\begin{document}

\maketitle


\begin{abstract}
We address Newton-type problems of minimal resistance
from an optimal control perspective.
It is proven that for Newton-type problems the Pontryagin maximum principle
is a necessary and sufficient condition.
Solutions are then computed for concrete situations,
including the new case when the flux of particles is non-parallel.
\end{abstract}


\smallskip


\smallskip


Pontryagin maximum principle, non-parallel flux of particles.


\section{Introduction}

In 1686, in his celebrated \emph{Principia Mathematica},
Isaac Newton propounded the problem of determining the profile
of a body of revolution, moving along its axis with constant speed,
through some resisting medium, which would minimize the
total resistance (see \cite{ButtazzoKawohl93,HabilitationLachand}).
Problems of this kind find application
in the building of high-speed and high-altitude flying vehicles, such as
in the design of missiles or artificial satellites.
Newton has given the correct answer to his problem, in
the situation of a ``rare'' medium of perfectly
elastic particles with constant mass and at equal distances from
each other, the resisting pressure at a surface point
of the body being proportional to the square of the normal component
of its velocity, but without explaining how he obtained it.
He didn't write, however, ``I have a great proof, but no space for it
in the margins of this book''. A proof ``from the Book'' was
waiting for the Pontryagin maximum principle.

When one writes the resistance force $\mathcal{R}$ associated to
Newton's problem,
\begin{equation*}
\mathcal{R}\left[\dot{x}(\cdot)\right]
= \int_0^T t \, \dfrac{1}{1+\dot{x}(t)^2} \,\mathrm{d}t \, ,
\end{equation*}
one obtains an integral functional of the type of those
studied throughout the history of the calculus of variations.
However, due to the restrictions on the derivatives of admissible trajectories,
$\dot{x}(t) \ge 0$, no satisfactory theory is available within the calculus
of variations framework (see \cite{alekseev,troutman,young}).
As first noticed by Legendre in 1788 (see \cite{belloniKawohl} and references therein),
without such restrictions on the derivatives the problem has no solution
(the infimum is zero), since one can obtain arbitrarily small values for
the integral resistance $\mathcal{R}\left[\dot{x}(\cdot)\right]$
by choosing a zig-zag function $x(\cdot)$ wildly oscillating,
with large derivatives in absolute value.
To make the problem physically consistent one must take into account
the monotonicity of the profile, and this means, as
was first remarked by V.M.~Tikhomirov (\textrm{cf.} \cite{alekseev,kvant}),
that Newton's problem belongs to optimal control:
\begin{gather}
\mathcal{R}\left[u(\cdot)\right] =
\int_0^T t \, \frac{1}{1+u(t)^2} \,  \mathrm{d}t \longrightarrow \min \, , \notag \\
\dot{x}(t) = u(t) \, , \quad   u(t) \ge 0 \, , \label{eq:NewtonsPrb}\\
x(0) = 0 \, , \quad x(T) = H \, . \notag
\end{gather}
Most part of the literature wrongly assume Newton's problem
to be ``one of the first applications of the calculus of variations''
but, in spite of this, the same literature correctly asserts
the birth of the calculus of variations:
1697, the publication date of the solution to the brachystochrone problem,
and not 1686, the publication date of the solution to Newton's problem
of minimal resistance.

In 1997 H. J. Sussmann and J. C. Willems,
in the beautiful paper \cite{300Sussmann},
defended the polemic thesis that the brachystochrone date 1697
marks not only the birth of the calculus of variations
but also the birth of optimal control. The truth seems to
be deeper: optimal control was born in 1686, before the
calculus of variations, with Newton's problem of minimal resistance.
The restriction on the control $u(\cdot)$,
which appear in Newton's problem \eqref{eq:NewtonsPrb},
is a common ingredient of the optimal control problems.
Such constraints appears naturally in practical engineering control
problems, and are treated with the Pontryagin maximum principle --
the central result of optimal control theory,
first conjectured by L.~S.~Pontryagin, and then
proved, in the late 1950's, by him and his collaborators,
V.~G.~Boltyanskii, R.~V.~Gamkrelidze, and E.~F.~Mishchenko
\cite{pontryagin}. In an optimal control problem,
the control functions take values on a set which is, in
general, not a vector space. This is precisely what happens
in Newton's classical problem \eqref{eq:NewtonsPrb},
and the reason why Newton's problem must be classified as an
optimal control problem, and not as a problem of the
calculus of variations.

Newton's problem has been widely studied,
and the literature about it is extensive.
The main difficulty is that of existence \cite{MR99b:49039}:
the Lagrangian $L(t,u) = \frac{t}{1+u^2}$ associated
to Newton's problem \eqref{eq:NewtonsPrb} is neither coercive nor convex,
and Tonelli's direct method (see \cite{MR85c:49001}) fails.
In order to prove existence, several different
classes of admissible functions have been proposed.
The question is now usually treated with the help of relaxation techniques
(see \cite{buttazzo99}), although direct arguments are also possible
(see \cite{plakhovCM03I05,plakhovCM03I07}).
As we shall prove (\S \ref{sec:MResult}), for the Newton-type problems,
the existence of a minimizer follows easily from the Pontryagin maximum principle:
one can show that the Pontryagin extremals are, for such problems,
absolute minimizers (\textrm{cf.} Theorem~\ref{res:PEequivAM}).

Several extensions of Newton's problem have been considered in
recent years. This revival of interest in Newton's problem,
and in the study of many variations around it, has
been motivated by the paper \cite{ButtazzoKawohl93}
of G.~Buttazzo and B.~Kawohl. Recent results
on Newton-type problems include:
bodies without rotational symmetry (nonsymmetric cases)
\cite{MR97i:49002,ButtazzoFeroneKawohl,lachand}; unbounded body
(resistance per unit area) with one-impact assumption
\cite{ComteAndLachandRobert}; bodies with rotational symmetry
and one-impact assumption, but not convex \cite{comte};
friction between particles and body (non-elastic collisions)
\cite{NoDEA2002}; bodies with prescribed volume \cite{belloniWagner};
multiple collisions allowed \cite{plakhovCM03I05,plakhovCM03I07};
unbounded body and multiple collisions allowed \cite{plakhovJDCS}.
More recently, Newton-type problems have been related with
problems of mass transportation \cite{plakhovCM03I34,plakhovCM03I36}.
For a good survey on mass optimization problems and open problems,
we refer the reader to \cite{MR2006304}.

Here we consider convex $d$-dimensional bodies of revolution
with Height $H$ and radius of maximal cross section $T$,
and treat them using an optimal control approach.
We will not be restricted to two-dimensional
or three-dimensional bodies,
considering bodies of arbitrary dimension $d \ge 2$.
We also introduce a different point of view. For us
the body does not move, and the particles are the ones who move.
The body is situated in a flux of infinitesimal particles,
the flux being invariant with respect to translations and rotations
around the symmetry axis of the body. This new point of view is,
in our opinion, physically more realistic. Newton has considered
the particles with no temperature (not moving). When the particles
have temperature, they move, and the flux of particles
is not necessarily falling vertically downwards the body,
as considered by Newton. We will be considering new interesting situations
with a non-parallel flux of particles. We obtain complete solution
to this class of Newton-type problems, by showing
that, under some physically relevant assumptions on the Lagrangian,
a control is an absolute minimizing control for the problem if, and only if,
it is a Pontryagin extremal control. Thus, for the Newton-type problems
we are dealing with, Pontryagin maximum principle holds
not only as a necessary optimality condition, but also
as a sufficient condition. As very special
situations, one obtains the solution found by Newton himself
(\S \ref{Sec:NewPrb}), and solutions to Newton's problem
in higher-dimensions (\S \ref{Sec:NewPrbHD}).


\section{Optimal Control}

The optimal control problem in Lagrange form consists in the minimization
of an integral functional
\begin{equation}
\label{eq:IntFunc}
J\left[x(\cdot),u(\cdot)\right] =
\int_0^T L\left(t,x(t),u(t)\right) \mathrm{d}t
\end{equation}
among all the solutions of a differential equation
\begin{equation}
\label{eq:controlEquation}
x'(t) = \varphi\left(t,x(t),u(t)\right) \, , \quad t \in [0,T] \,
\end{equation}
subject to the boundary conditions
\begin{equation}
\label{eq:boundaryCond}
x(0) = \alpha \, , \quad x(T) = \beta \, .
\end{equation}
The Lagrangian $L$ and the velocity function $\varphi$ are
defined on $[a,b] \times \mathbb{R}^n \times \Omega$,
where $\Omega \subseteq \mathbb{R}^r$ is called the control set.
The main difference between the problems of optimal control
and those of the calculus of variations, is that $\Omega$
is in general not an open set. In the case $\varphi(t,x,u) = u$,
and $\Omega = \mathbb{R}^n$, one gets the fundamental problem of the
calculus of variations. For the Newton problem, we have $n = r = 1$,
$\Omega = \mathbb{R}_{0}^{+}$, $\varphi(t,x,u) = u$, $\alpha = 0$,
$\beta > 0$, and $L(t,x,u) = \frac{t}{1+u^2}$.
Typically, $L(t,x,u)$ and $\varphi(t,x,u)$
are continuous with respect to all arguments and have continuous
derivatives with respect to $x$; the admissible processes
$\left(x(\cdot),u(\cdot)\right)$ are formed by absolutely continuous
state trajectories $x(\cdot)$ and measurable and bounded controls
$u(\cdot)$, taking values on the control set $\Omega$ and
satisfying \eqref{eq:controlEquation}-\eqref{eq:boundaryCond}.

The Pontryagin maximum principle is a first-order necessary optimality
condition, which provides a generalization of the classical Euler-Lagrange
equations and Weierstrass condition, to problems in which upper and/or lower
bounds are imposed on the control variables.

\begin{theorem}[Pontryagin maximum principle]
\label{th:PMP}
Let $\left(x(\cdot),u(\cdot)\right)$ be a minimizer of the optimal
control problem. Then there exists a pair
$\left(\psi_0,\psi(\cdot)\right)$,
where $\psi_0 \le 0$ is a constant
and $\psi(\cdot)$ an $n$-vector absolutely continuous function with
domain $[0,T]$, \emph{not all zero},
such that the following holds true for almost all $t$ on
the interval $[0,T]$:
\begin{description}
\item[(i)] the Hamiltonian system
\begin{equation*}
\begin{cases}
x'(t) =
\frac{\partial \mathcal{H}}{\partial \psi}\left(t,x(t),u(t),\psi_0,\psi(t)\right) \, , \\
\psi'(t) =
- \frac{\partial \mathcal{H}}{\partial x}\left(t,x(t),u(t),\psi_0,\psi(t)\right) \, ;
\end{cases}
\end{equation*}
\item[(ii)] the maximality condition
\begin{equation}
\label{PMP:MC}
\mathcal{H}\left(t,x(t),u(t),\psi_0,\psi(t)\right) = \max_{u \in \Omega}
\mathcal{H}\left(t,x(t),u,\psi_0,\psi(t)\right) \, ;
\end{equation}
\end{description}
where the Hamiltonian $\mathcal{H}$ is defined by
\begin{equation*}
\mathcal{H}(t,x,u,\psi_0,\psi) = \psi_0 L(t,x,u) + \psi \cdot \varphi(t,x,u) \, .
\end{equation*}
\end{theorem}

The first equation in the Hamiltonian system is just the control equation
\eqref{eq:controlEquation}. The second equation is known as the
\emph{adjoint system}.

\begin{definition}
A quadruple $\left(x(\cdot),u(\cdot),\psi_{0},\psi(\cdot)\right)$
satisfying the Hamiltonian system and the maximality condition
is called a Pontryagin extremal. The control $u(\cdot)$ is said to be
an extremal control. The extremals are said to be
abnormal when $\psi_0 = 0$ and normal otherwise.
\end{definition}

\begin{remark}
If $\left(x(\cdot),u(\cdot),\psi_{0},\psi(\cdot)\right)$
is a Pontryagin extremal, then, for any $\gamma > 0$,
$\left(x(\cdot),u(\cdot),\gamma \psi_{0},\gamma \psi(\cdot)\right)$
is also a Pontryagin extremal. From this simple observation
one can consider, without any loss of generality, that
$\psi_0 = -1$ in the normal case.
\end{remark}

\begin{remark}
The fact that Theorem~\ref{th:PMP} asserts the existence
of Hamiltonian multipliers $\psi_0$ and $\psi(\cdot)$
not vanishing simultaneously is of primordial
importance: without this condition, all admissible pairs
$\left(x(\cdot),u(\cdot)\right)$ would be Pontryagin
extremals.
\end{remark}

In some situations, it may happen that functions $L$ and/or
$\varphi$ depend upon some parameters $w \in W \subseteq
\mathbb{R}^k$. In this case, given a control $u(\cdot)$, the
corresponding state trajectory $x(\cdot)$ and the cost functional
$J$ will in general depend on the choice of the parameters
$w$. The problem in then to choose the parameters $\tilde{w}$ in
$W$ for which there exists an admissible pair
$\left(\tilde{x}(\cdot),\tilde{u}(\cdot)\right)$ such that
$J\left[\tilde{x}(\cdot),\tilde{u}(\cdot),\tilde{w}\right] \le
J\left[x(\cdot),u(\cdot),w\right]$ for all $w \in W$ and
corresponding admissible pairs $\left(x(\cdot),u(\cdot)\right)$.
The parameter problem can be reformulated in the format
\eqref{eq:IntFunc}--\eqref{eq:controlEquation} by considering $w$ as
a state variable with dynamics $w'(t) = 0$ and initial
condition $w(0) \in W$.


\section{Optimal Control of Newton-Type Problems}
\label{sec:MResult}

The standard method to solve a problem in optimal control
proceeds by first proving that a solution to the problem exists;
then assuring the applicability of the Pontryagin maximum principle;
and, finally, identifying the Pontryagin extremals (the candidates).
Further elimination, if necessary, identifies the minimizer or
minimizers of the problem. It is not easy to prove existence
for Newton's problem with the classical arguments, because the
Lagrangian $L(t,u) = \frac{t}{1+u^2}$ is not coercive and
it is not convex with respect to $u$ for $u \ge 0$. Here we
will make use of a different approach. We will show, by a simple
and direct argument, that for Newton-type problems
\eqref{eq:HypL}--\eqref{NTP} the Pontryagin extremals are
absolute minimizers. This means that, in order
to solve a Newton-type problem, it is enough to identify
the Pontryagin extremals (\textrm{cf.} Theorem~\ref{res:PEequivAM}).

We begin to show that there are no abnormal extremals
for a Newton-type problem.

\begin{proposition}
\label{Prop:NoAbnormality}
Let $L(t,u)$ be a continuous function satisfying the following
conditions:
\begin{gather}
L(t,u) > \xi \ge 0 \quad \forall \, (t,u) \in \,
]0,T] \times \mathbb{R}_{0}^{+} \, , \label{eq:HypL} \\
\lim_{u \rightarrow +\infty} L(t,u) = \xi  \quad \forall \, t \in [0,T] \, . \notag
\end{gather}
Then all Pontryagin extremals $(x(\cdot),u(\cdot),\psi_0,\psi(\cdot))$
of the Newton-type problem
\begin{gather}
J\left[u(\cdot)\right] =
\int_0^T L\left(t,u(t)\right) \mathrm{d}t \longrightarrow \min \, , \notag \\
x'(t) = u(t) \, , \quad u(t) \ge 0 \, , \label{NTP} \\
x(0) = 0 \, , x(T) = \beta \text{ with } \beta > 0 \, , \notag
\end{gather}
are normal extremals ($\psi_0 = -1$) with $\psi$ a negative constant
($\psi(t) \equiv - \lambda$, $\lambda > 0$).
\end{proposition}

\begin{proof}
As far as the Hamiltonian does not depend on $x$,
\begin{equation*}
\mathcal{H}\left(t,u,\psi_0,\psi\right) = \psi_0 L(t,u) + \psi u \, ,
\end{equation*}
we conclude from the adjoint system that $\psi(t) \equiv c$,
with $c$ a constant. If $c$ is equal to zero, then $\psi_0 < 0$
(they are not allowed to be both zero) and the maximality
condition \eqref{PMP:MC} simplifies to
\begin{equation*}
\psi_0 L(t,u(t)) = \max_{u \ge 0}
\left\{\psi_0 L(t,u)\right\} \, .
\end{equation*}
Under the hypotheses \eqref{eq:HypL} the maximum is not
achieved ($u \rightarrow +\infty$) and we conclude that
$c \ne 0$. Similarly, for $c > 0$ the maximum
\begin{equation*}
\max_{u \ge 0} \left\{\psi_0 L(t,u) + c u\right\}
\end{equation*}
does not exist and one concludes that $c < 0$. It remains
to prove that $\psi_0$ is different from zero. Indeed, if
$\psi_0 = 0$, \eqref{PMP:MC} reads
\begin{equation*}
c u(t) = \max_{u \ge 0} \left\{c u \right\} \, ,
\end{equation*}
and follows that $u(t) \equiv 0$ and
$x(t) \equiv \text{constant}$. This is not a possibility
since $\beta > 0$. The proof is complete.
\end{proof}

Theorem~\ref{res:PEequivAM} reduces the
procedure of solving a Newton-type problem to the
computation of Pontryagin extremals.

\begin{theorem}
\label{res:PEequivAM}
The control $\tilde{u}(\cdot)$ is an absolute minimizing control
for the Newton-type problem \eqref{eq:HypL}--\eqref{NTP},
\textrm{i.e.}, $J[\tilde{u}(\cdot)] \le J[u(\cdot)]$ for all
$u(\cdot) \in L_{\infty}\left([0,T],\mathbb{R}_{0}^{+}\right)$,
if, and only if, it is an extremal control.
\end{theorem}

\begin{proof}
Theorem~\ref{res:PEequivAM} is a direct consequence of the maximality
condition. From Proposition~\ref{Prop:NoAbnormality} one can write
\eqref{PMP:MC} as
\begin{equation}
\label{eq:FromMC}
- L(t,\tilde{u}(t)) - \lambda \tilde{u}(t) \ge
- L(t,u(t)) - \lambda u(t)
\end{equation}
for all admissible controls $u(\cdot)$
and for almost all $t \in [0,T]$.
Having in mind that all admissible processes $\left(x(\cdot),u(\cdot)\right)$
of the problem \eqref{NTP} satisfy
\begin{equation*}
\int_0^T u(t) \mathrm{d}t = \int_0^T x'(t) \mathrm{d}t = \beta \, ,
\end{equation*}
it is enough to integrate \eqref{eq:FromMC} to obtain the
desirable conclusion:
\begin{equation*}
\int_0^T L(t,\tilde{u}(t)) \mathrm{d}t
\le \int_0^T L(t,u(t)) \mathrm{d}t \, .
\end{equation*}
\end{proof}

The required optimal solutions of the Newton-type problem
\eqref{eq:HypL}--\eqref{NTP} are exactly the
Pontryagin extremals. This means, essentially, that we have
reduced a dynamic optimization problem (a minimization problem
in the space of functions) to the static optimization problem
given by the maximality condition.

\begin{corollary}
\label{cor:DOtoSO}
Finding the solutions for the Newton-type problem
\eqref{eq:HypL}--\eqref{NTP} amounts to find the minimum of
the function $h(u) = L(t,u) + \lambda u$, $t \in [0,T]$,
$\lambda > 0$, for $u \ge 0$.
\end{corollary}


\section{An Application}
\label{sec:App}

Consider a $d$-dimensional body of revolution
$$
\{ (\xi_0,\, \xi): \xi_0 \in [0,\, H],\ |\xi| < \Phi(\xi_0) \}
\subset {\mathbb R}^{d},
$$
where $d \ge 2$,\, $\xi = (\xi_1, \ldots, \xi_{d-1})$,\, $\Phi$ is
a non-negative function defined on $[0,\, H]$. Denote by $T$ the
radius of maximal cross section of the body, $T =
\max_{0\le\zeta\le H} \Phi(\zeta)$. Let us assume that the body is
convex, then the function $\Phi$ is concave, and there exists $c
\in [0,\, H]$ such that $\Phi(\zeta)$ is monotone increasing as
$\zeta \le c$, and monotone decreasing as $\zeta \ge c$.

We suppose that the body is unmovable and is situated in a flux of
infinitesimal particles. The flux is invariant with respect to
translations and rotations around the $\xi_0$-axis, which is the
symmetry axis of the body. So, the specific pressure of the flux
on an infinitesimal element of the body surface depends only on
the value of $\Phi'$ at that element. It is convenient to
consider, instead of $\Phi$, two functions that are generalized
inverses of $\Phi$; denote them by $x_-(t)$ and $H - x_+(t)$. They
are defined in the following way:\ $x_-(t) = 0$ as $t \in [0,\,
\Phi(0)]$,\, and $x_-(t)$ is inverse to the strictly monotone
increasing branch of $\Phi$ as $t \in [\Phi(0),\, T]$;\ $x_+(t) =
0$ as $t \in [0,\, \Phi(H)]$,\, and $H - x_+(t)$ is inverse to the
strictly monotone decreasing branch of $\Phi$ as $t \in
[\Phi(H),\, T]$. The obtained functions $x_-$ and $x_+$ are
convex, continuous, and monotone increasing, besides $x_-(0) =
x_+(0) = 0$,\, $x_-(T) \le c$,\, $x_+(T) \le H - c$. In such a
representation, the specific pressure is a function of $x_+'$ or
of $x_-'$, if the point belongs to the front or to the rear part
of surface, respectively; we denote the corresponding functions by
$p_+(\cdot)$ and $-p_-(\cdot)$.

The pressure on an element $d^{d-1} s$ of the front part of
surface is $dp_+ = p_+(x_+') d^{d-1} s$. The projection of the
pressure vector to the $\xi_0$-axis equals $dp_0 = dp_+/ \sqrt{1 +
x_+'^2}$, and the projection of the surface element to
${\mathbb R}^{d-1}_{\xi_1\ldots \xi_{d-1}}$ has area $d^{d-1} \xi =
d^{d-1} s/ \sqrt{1 + x_+'^2}$. Thus, the $\xi_0$-projection of
pressure corresponding to the element $d^{d-1} \xi$ is $dp_0 =
p_+(x_+')d^{d-1} \xi$. Passing to polar coordinates and
integrating over the ball $\{ |\xi| < T \}$, one obtains the
resistance ${\cal R}_+$ of the front part of body to the flux:
$$
{\cal R}_+[x_+(\cdot)] = v_{d-1} \int_0^T\, p_+(x_+'(t))\, dt^{d-1},
$$
here $v_{d-1}$ stands for the volume of $(d-1)$-dimensional unit
ball. Similarly, the resistance of the rear part of body to the
flux (which is positive) equals $-{\cal R}_-[x_-(\cdot)]$, where
$$
{\cal R}_-[x_-(\cdot)] = v_{d-1} \int_0^T\, p_-(x_-'(t))\, d t^{d-1}.
$$
So, the resistance of body to the flux is ${\cal R} [x_+(\cdot),
x_-(\cdot)] = {\cal R}_+[x_+(\cdot)] + {\cal R}_-[x_-(\cdot)]$.

It is required to minimize ${\cal R} [x_+(\cdot), x_-(\cdot)]$ over
all pairs $(x_+(\cdot),\, x_-(\cdot))$ of convex monotone
increasing functions defined on $[0,\, T]$,\, provided $x_\pm$
take values in $[0,\, \beta_\pm]$, where $\beta_- = c$,\, $\beta_+ = H -
c$,\, $T$ and $H$ are fixed, and $c$ varies between 0 and $H$.

We are acting as follows. First we fix the sign "+" or "$-$",
minimize ${\cal R}_\pm$ over monotone increasing functions $x: [0,\,
T] \mapsto [0,\, \beta_\pm]$,\, with $\beta_\pm$ fixed, and verify
that among all the solutions, the convex one is unique; denote it
by $x_\pm^{\beta_\pm}$. Then we minimize the sum
${\cal R}_+[x_+^{\beta_+}(\cdot)] + {\cal R}_-[x_-^{\beta_-}(\cdot)]$ over
all positive $\beta^+$ and $\beta^-$ such that $\beta^+ + \beta^- = H$.


\subsection{Solving the problem in general case}
\label{Sec:GenCase}

In what follows, we assume that the functions $p_+$ and $p_-$
satisfy the following conditions:

(i) $p_\pm \in C^1 [0,\, +\infty)$;

(ii) there exist $\lim_{u\to +\infty} p_\pm(u)$;

(iii) $p_\pm'(0) = \lim_{u\to +\infty} p_\pm'(u) = 0$;

(iv) for some $\bar{u}_\pm > 0$,\, $p_\pm'$ is strictly monotone
decreasing on $[0,\, \bar{u}_\pm]$, and strictly monotone
increasing on $[\bar{u}_\pm,\, +\infty)$. \vspace{1mm}

For simplicity, we further put $v_{d-1} = 1$. Let us fix the sign
"+" or "$-$", and introduce shorthand notations
$$
p_\pm = p, \ \ \ \beta_\pm = \beta, \ \ \ {\cal R}_\pm = {\cal R}, \ \ \
x_\pm^{\beta_\pm} = x^\beta.
$$

\begin{proposition}\label{sol}
There exists a unique solution $u^0$ of the problem
$$
\frac{p(0) - p(u)}{u} \to \max,
$$
besides $u^0 > \bar{u}$.
\end{proposition}

\begin{proof}
Denote $q(u) = p(0) - p(u)$ and $B = \sup_{u>0} (q(u)/ u)$. From
(i)--(iv) it follows that the function $q(u)/ u$,\, $u
> 0$\, is continuous, positive, and satisfies the relations $\lim_{u\to 0+}
(q(u)/ u) = \lim_{u\to +\infty} (q(u)/ u) = 0$, hence $0 < B <
+\infty$, and there exists a value $u^0 > 0$ such that $q(u^0)/
u^0 = B$. Obviously, at $u = u^0$ one has $\left( q(u)/ u \right)'
= 0$, hence $q'(u^0) = q(u^0)/ u^0$. At some $\theta \in (0,\, 1)$
one has $q(u^0)/ u^0 = q'(\theta u^0)$, hence $q'(u^0) = q'(\theta
u^0)$. This implies that $q'$ is not strictly monotone on $[0,\,
u^0]$;\, thus, by virtue of (iv), $u^0 > \bar{u}$.

It remains to prove that the value $u^0$, solving the equation
$q(u)/ u = B$, is unique. Suppose that $q(u^0)/ u^0 = q(u^1)/ u^1
= B$,\, $u^0 < u^1$.\, Then $q(u^0) = u^0\, q'(u^0)$,\, $q(u^1) =
u^1\, q'(u^1)$. At some $u \in (u^0,\, u^1)$,\, one has $q(u^1) -
q(u^0) = q'(u)\, (u^1 - u^0)$;\, this implies that
$$
q'(u)\, (u^1 - u^0) = u^1\, q'(u^1) - u^0\, q'(u^0),
$$
hence
\begin{equation}\label{positiv}
u^0\, (q'(u^0) - q'(u)) + u^1\, (q'(u) - q'(u^1)) = 0.
\end{equation}
One has $u^0 > \bar{u}$, hence $q'$ is strictly monotone
decreasing as $u \ge u^0$, so both terms in (\ref{positiv}) are
positive. The obtained contradiction proves the proposition.
\end{proof}

Let us denote
$$
B = \frac{p(0) - p(u)}{u} = - p'(u).
$$

\begin{proposition}\label{sol2}
(a) As $\lambda\, t^{2-d} > B$, the unique solution of the problem
\begin{equation}
\label{condition} t^{d-2} \, p(u) + \lambda u \to \min;
\end{equation}
is $u = 0$.

(b) As $\lambda\, t^{2-d} = B$, there are two solutions:\, $u = 0$
and $u = u^0$.

(c) As $\lambda\, t^{2-d} < B$, the solution $\tilde u$ is unique,
besides $\tilde u > u^0$, and $p'(\tilde u) = -\lambda t^{2-d}$.
\end{proposition}

\begin{proof}
(a) and (b) are obvious; let us prove (c). Denote $\tilde \lambda :=
\lambda\, t^{2-d}$. By definition of $B$,\, for $0 < u < u^0$ one has
\begin{gather*}
\frac{p(0) - p(u)}{u} < B = \frac{p(0) - p(u^0)}{u^0} \, , \\
p(u^0) + B u^0 = p(0) < p(u) + B u \, ,
\end{gather*}
hence
$$
p(u) - p(u^0) > B\, (u^0 - u) > \tilde \lambda\, (u^0 - u),
$$
and thus,
$$
p(u) + \tilde \lambda u > p(u^0) + \lambda u^0.
$$
On the other hand, one has $B = -p'(u^0)$, therefore
$$
p'(u^0) + \tilde \lambda < 0;
$$
moreover the function $p(u) + \tilde \lambda u$ is convex on $[u^0,\,
+\infty)$ and tends to infinity as $u \to +\infty$. All this
implies that the solution $\tilde u$ of (\ref{condition}) is
unique, satisfies the equation $p'(\tilde u) + \tilde \lambda = 0$,
and $\tilde{u} > u^0$.
\end{proof}

From Corollary~\ref{cor:DOtoSO} we know that if $x^\beta(\cdot)$ is
a solution of the minimization problem $\, {\cal R}[x(\cdot)] \to
\min$,\, $x: [0,\, T] \to [0,\, \beta]$,\, then for some $\lambda$,\,
the values $u = {x^{\beta\,}}'(t)$,\, $t \in [0,\, T]$ satisfy the
equation (\ref{condition}). According to propositions \ref{sol}
and \ref{sol2}, one should distinguish between three cases:\,
(a)\, if $\lambda t^{2-d} > B$, then $u = 0$;\, (b)\, if $\lambda
t^{2-d} = B$, then $u = 0$ or $u = u^0$;\, (c)\, if $\lambda t^{2-d}
< B$, then $u > u^0$, and $p'(u) = -\lambda t^{2-d}$.

Consider two different cases: $d = 2$ (two-dimensional problem)
and $d \ge 3$ (the problem in three or more dimensions).

\subsubsection{Two-dimensional problem ($d = 2$)}
\label{d=2}

If $\lambda > B$, the unique solution of (\ref{condition}) is $u =
0$, hence $x^\beta \equiv 0$. This implies that $\beta = 0$.

If $\lambda = B$, there are two solutions:\, $u = 0$ and $u = u^0$,
therefore any absolutely continuous function $x(\cdot)$,\, $x(0) =
0$,\, $x(T) = \beta$,\, such that $x'(t)$ takes the values 0 and
$u^0$, minimizes ${\cal R}$. A convex solution $x^\beta$ has monotone
increasing derivative, hence for some $t_0 \in [0,\ T]$,\,
${x^{\beta\,}}'(t) = 0$ as $t \in [0,\ t_0]$,\, and ${x^{\beta\,}}'(t)
= u^0$ as $t \in [t_0,\, T]$.\, Thus,
\begin{equation}\label{xt}
x^\beta(t) = \left\{ \begin{array}{lc} 0 & \text{ as } t \in [0,\ t_0]\\
u^0\, (t - t_0)  & \text{ as } t \in [t_0,\, T]. \end{array}
\right.
\end{equation}
Taking into account that $x^\beta(T) = \beta$, one concludes that
$\beta/ T \le u^0$ and $t_0 = T - \beta/ u^0$.

If $\lambda < B$, there is a unique solution $\tilde u$,\, hence
$x^\beta(t) = \tilde u t$,\, and $\beta/ T = \tilde u > u^0$.
  \vspace{2mm}

Summarizing, one gets that\\
(i)\ as $\beta/ T < u^0$,\, the convex solution $x^\beta(t)$ is given by (\ref{xt});\\
(ii)\ as $\beta/ T \ge u^0$,\,\ $x^\beta(t) = \beta\, t/ T$.
  \vspace{2mm}

As $\beta/ T < u^0$, one has
$$
{\cal R}[x^\beta(\cdot)] = \int_0^{t_0} p(0)\, dt + \int_{t_0}^T
p(u^0)\, dt,
$$
and taking into account that $t_0 = T - \beta/ u^0$ and $(p(0) -
p(u^0)) /u^0 = B$, one gets
$$
{\cal R}[x^\beta(\cdot)] = T\, p(0) - \beta\, B.
$$
As $\beta/ T \ge u^0$,\, one has ${\cal R}[x^\beta(\cdot)] = T\,
p(\beta/ T)$. Introduce the function
$$
\bar{p}(u) = \left\{ \begin{array}{ll} p(0) - B\, u, & \ \ \text{ if }
 \ 0 \le u \le u^0,\\
p(u), & \ \ \text{ if } \ u \ge u^0,
\end{array} \right.
$$
then
$$
{\cal R}[x^\beta(\cdot)] = T\, \bar{p}(\beta/ T).
$$
Thus, the minimization problem ${\cal R}_+[x_+^{\beta_+}(\cdot)] +
{\cal R}_-[x_-^{\beta_-}(\cdot)] \to \min$ is reduced to the problem
\begin{equation}\label{problem3passo}
p_h(z) = \bar{p}_+(z) + \bar{p}_-(h - z) \to \min, \ \ \ \ \ \ 0
\le z \le h,
\end{equation}
where $h = H/T$.
 \vspace{2mm}

The introduced functions $\bar p_\pm(u)$ are continuously
differentiable on $[0,\, +\infty)$, and
$$
\bar{p}_\pm'(u) = \left\{ \begin{array}
{ll} -B_\pm & \ \ \text{ if } \ 0 \le u \le u_\pm^0,\\
p'_\pm(u) & \ \ \text{ if } \ u > u_\pm^0.
\end{array} \right.
$$
Using that $u_\pm^0 > \bar u_\pm$, one concludes that $\bar
p_\pm'(u)$ is monotone increasing, hence $p_h'(z)$,\, $0 \le z \le
h$\, is also monotone increasing.

From now and until the end of subsection \ref{Sec:GenCase}, we
shall assume that \label{assump:pg:pMq} $p'_+(u) < p'_-(u)$,\, $u
\ge 0$,\, hence $B_+ > B_-$.\, Denote by $u_*$ a positive value
such that $\bar p_+'(u_*) = -B_-$. This value is unique, and $u_*
> u_+^0$.\,
Consider four cases:\\
1) $0 < h < u_+^0$; \ 2) $u_+^0 \le h \le u_*$; \ 3) $u_* < h <
u_* + u_-^0$; \ and 4) $h \ge u_* + u_-^0$.

In the cases 1) and 2) one has
$$
p_h'(z) \le p_h'(h) = \bar p_+'(h) + B_- \le 0 \ \ \ \text{ as } \
\ 0 \le z \le h,
$$
hence the minimum of (\ref{problem3passo}) is achieved at $z = h$.
Therefore, the optimal value of $\beta_-$ is zero, so $x_-^{\beta_-}
\equiv 0$.

In the case 1)\, one has $\beta_+ /T = h < u_+^0$, hence
$x_+^{\beta_+}(\cdot)$ is given by (\ref{xt}), with $t_0 = T\, (1 -
h/u_+^0)$.\, So, the optimal body is a trapezium.

In the case 2)\, one has $x_+^{\beta_+}(t) = ht$, hence the optimal
body is an isosceles triangle.

In the cases 3) and 4),\, one has $\bar p_+'(h) > -B_- > -B_+ =
\bar p_+'(u_+^0)$, hence $h > u_+^0$.\, Further, one has
$$
p_h'(h) = \bar p_+'(h) - B_- > 0;
$$
on the other hand,
$$
p_h'(u_+^0) = \bar p_+'(u_+^0) - \bar p_-'(h - u_+^0) \le -B_+ +
B_- < 0.
$$
It follows that the minimum of $p_h$ is achieved at an interior
point of $[u_+^0,\, h]$, so the optimal value of $\beta_-$ satisfies
the relation $u_+^0 < \beta_+/T < h$,\, and $x_+^{\beta_+}(t) = t\,
\beta_+/T$.

In the case 3),\, denoting $\tilde h = \max \{ 0,\, h - u_-^0 \}$,
one has $\tilde h < u_*$, hence
$$
p_h'(\tilde h) = \bar p_+'(\tilde h) - \bar p_-'(h - \tilde h) \le
\bar p_+'(\tilde h) + B_- < 0,
$$
hence the minimum of $p_h$ is reached at an interior point of
$[\tilde h,\, h]$,\, thus $0 < \beta_-/T < h - \tilde h \le u_-$,
and
$$
x_-^{\beta_-}(t) = \left\{ \begin{array}{ll} 0 & \text{ if } t \in [0,\ T - \beta_-/u_-^0]\\
u_-^0\, (t - T + \beta_-/u_-^0) & \text{ if } t \in [T -
\beta_-/u_-^0,\, T].
\end{array} \right.
$$
The optimal body here is the union of a triangle and a trapezium
turned over.

In the case 4),\, one has $p_h'(h - u_-^0) = \bar p_+'(h - u_-) +
B_- \ge 0$,\, hence the minimum of $p_h$ is reached at a point of
$[u_+^0,\, h - u_-^0)$.\, Thus, $\beta_-/T > u_-^0$, and
$x_-^{\beta_-}(t) = t\, \beta_-/T$.\, The optimal body is a union of
two isosceles triangles with common base.

\subsubsection{Problem in three or more dimensions ($d \ge 3$)}
\label{d ge 3}

Here we additionally assume that $p_\pm \in C^2 [0,\, +\infty)$
and $p_+''(u) > 0$ as $u > u_+^0$.

Denote $\omega = 1/(d - 2)$ and $t_0 = (\lambda/B)^\omega$. As $0 \le t <
t_0$, the unique solution of (\ref{condition}) is $u = 0$, hence
$x^\beta(t) = 0$.\, As $t_0 < t \le T$, the solution $u$ satisfies
the relation
$$
t^{d-2}\, p'(u) + \lambda = 0,
$$
and $u > \tilde u$.

If $T \le t_0$,\, one has $x^\beta \equiv 0$ and $\beta = 0$.\, Let,
now, $t_0 < T$;\, using that $p'$ is negative, continuous, and
strictly monotone increasing on $[u^0,\, +\infty)$, one concludes
that ${x^{\beta\,}}'(t) = u$ is also continuous, and is strictly
monotone increasing on $[t_0,\, T)]$ from ${x^{\beta\,}}'(t_0+) =
u^0$ to the value $U$ defined from the relation
\begin{equation*}
T^{d-2}\, p'(U) + \lambda = 0, \ \ \ \ U > u^0.
\end{equation*}
Thus, $x^\beta(\cdot)$ is convex;\, moreover, as $t_0 \le t \le
T$,\, $x^\beta$ can be represented as a function of $u \in [u^0,\,
U]$. Using that ${x^{\beta\,}}'(t) = u$ and that
\begin{equation}\label{111}
t = \frac{\lambda^\omega}{|p'(u)|^\omega},
\end{equation}
one gets
$$
\frac{dx^\beta}{du} = \frac{dx^\beta}{dt} \frac{dt}{du} = u\,
\lambda^\omega\, \frac{d}{du} \left( \frac{1}{|p'(u)|^\omega} \right),
$$
hence
$$
x^\beta = \lambda^\omega\, \int_{u^0}^u\, \nu\, d \left(
\frac{1}{|p'(\nu)|^\omega} \right);
$$
using that $|p'(u^0)| = B$, one obtains
\begin{equation*}
x^\beta = \lambda^\omega \left( \frac{u}{|p'(u)|^\omega} - \frac{u^0}{B^\omega}
- \int_{u^0}^u \frac{d\nu}{|p'(\nu)|^\omega} \right).
\end{equation*}
In particular, substituting $u = U$, one has
\begin{equation}\label{btbt}
\lambda^\omega\, \left( \frac{U}{|p'(U)|^\omega} - \frac{u^0}{B^\omega} -
\int_{u^0}^U \frac{d\nu}{|p'(\nu)|^\omega} \right) = \beta.
\end{equation}
Introduce the function
$$
g(u) = \int_0^u\, \frac{d\nu}{|\bar p'(\nu)|^\omega}.
$$
Using that $|\bar p'(\nu)| = B$ as $0 \le \nu \le u^0$, and $\bar
p'(\nu) = p'(\nu)$ as $\nu \ge u^0$, one gets
$$
g(U) = \frac{u^0}{B^\omega} + \int_{u^0}^U
\frac{d\nu}{|p'(\nu)|^\omega},
$$
and using that
$$
T = \frac{\lambda^\omega}{|p'(U)|^\omega},
$$
from (\ref{btbt}) one gets
\begin{equation}\label{btbtbt}
\frac{\beta}{T} = U - |p'(U)|^\omega\, g(U).
\end{equation}

The minimal resistance equals
\begin{equation*}
{\cal R}[x^\beta(\cdot)] = \int_0^T\, p(u(t))\, dt^{d-1} = p(0)\,
t_0^{d-1} + \int_{t_0}^T\, p(u(t))\, dt^{d-1}.
\end{equation*}
Using that $u(T) = U$,\, $u(t_0) = u^0$,\, $|p'(u^0)| = B$,\,
$p(0) - p(u^0) = B\, u^0$,\, $t_0 = \lambda^\omega/ B^\omega$, and also the
formula (\ref{111}), one finds
\begin{equation*}
\begin{split}
{\cal R}[x^\beta(\cdot)] &= \lambda^{1+\omega}\, \left\{
\frac{p(0)}{B^{1+\omega}} + \frac{p(U)}{|p'(U)|^{1+\omega}} -
\frac{p(u^0)}{B^{1+\omega}} - \int_{u^0}^U\,
\frac{dp(u)}{|p'(u)|^{1+\omega}} \right\} \\
&= \lambda^{1+\omega}\, \left\{ \frac{u^0}{B^\omega} +
\frac{p(U)}{|p'(U)|^{1+\omega}} + \int_{u^0}^U\,
\frac{du}{|p'(u)|^\omega} \right\}.
\end{split}
\end{equation*}
This implies
\begin{equation}\label{rrr}
\frac{{\cal R}[x^\beta(\cdot)]}{T^{d-1}} = p(U) + |p'(U)|^{1+\omega}\,
g(U).
\end{equation}
 \vspace{2mm}

Denote $U_+ = z_+$,\, $U_- = z_-$. Using (\ref{btbtbt}) and
(\ref{rrr}), one comes to the following problem of conditional
minimum
\begin{multline*}
r(z_-, z_+) := \\
p_+(z_+) + |p_+'(z_+)|^{1+\omega}\, g_+(z_+)\,   +\,   p_-(z_-) +
|p_-'(z_-)|^{1+\omega}\, g_-(z_-) \to \min,
\end{multline*}
under the conditions
\begin{equation}\label{conditiondim3}
z_- - |p_-'(z_-)|^\omega\, g_-(z_-)\, +\, z_+ - |p_-'(z_+)|^\omega\,
g_+(z_+) = h, \ \ \ z_- \ge u^0, \ \ z_+ \ge u^0.
\end{equation}
From (\ref{conditiondim3}), taking into account that
$|p_\pm'(z_\pm)|^\omega\, g_\pm'(z_\pm) = 1$, one obtains that $z_+$
is a differentiable function of $z_-$, and
$$
\frac{dz_+}{dz_-} = -\frac{|p_-'(z_-)|^{\omega-1}\, p_-''(z_-)\,
g_-(z_-)}{|p_+'(z_+)|^{\omega-1}\, p_+''(z_+)\, g_+(z_+)}.
$$
Now,
\begin{equation*}
\begin{split}
\frac{d}{dz_-} & \, r(z_-, z_+(z_-)) = \\
&= -\frac{dz_+}{dz_-}
\cdot (1 + \omega)\, |p_+'(z_+)|^\omega\, p_+''(z_+)\, g_+(z_+)\,  -\,
(1 + \omega)\, |p_-'(z_-)|^\omega\, p_-''(z_-)\, g_-(z_-)\\
&= (1 + \omega)\, |p_-'(z_-)|^{\omega-1}\, p_-''(z_-)\, g_-(z_-) \cdot
(p_-'(z_-) - p_+'(z_+)).
\end{split}
\end{equation*}
Note that $z_+$ is a monotone decreasing function of $z_-$, hence
the function $p_-'(z_-) - p_+'(z_+(z_-))$ is monotone increasing
as $z_- \ge u^0$,\, $z_+(z_-) \ge u^0$.

Recall that $u_*$ is the value satisfying $\bar p_+'(u_*) = -B_-$.
Denote
$$
h_* = u_* - B_-^\omega\, g_+(u_*).
$$
Consider two cases.

1)\ $h \le h_*$.\ One has $z_+(u^0) \le u_*$, hence $p_-'(u^0) -
p_+'(z_+(u^0)) \ge -B_- - p_+'(u_*) = 0$. It follows that as $z_-
> u^0$,\ $p_-'(z_-) - p_+'(z_+(z_-)) > 0$,\ so the minimum of $r(z_-,
z_+)$is attained at $z_- = 0$.

2)\ $h > h_*$.\ One has $z_+(u^0) > u_*$, hence $p_-'(u^0) -
p_+'(z_+(u^0)) < 0$.\, On the other hand, as $\tilde z_- =
z_+(\tilde z_-)$, one has $p_-'(\tilde z_-) - p_+'(\tilde z_-) <
0$, hence at some $z_- \in (u^0,\, \tilde z_-)$,\ $p_-'(z_-) =
p_+'(z_+(z_-))$,\, and so, the minimum of resistance is attained.


\subsection{Examples}

We have given in \S \ref{Sec:GenCase}
complete description of the solutions
to the formulated Newton-type problem.
We now consider, for illustration purposes,
various particular cases of the problem.
All the calculations can easily be done
with the help of a computer algebra system.
We have used Maple to implement a procedure
which, given functions $p_+(\cdot)$ and $p_-(\cdot)$
and the values for $T$ and $H$, gives the optimal shape
for the respective problem.

\subsubsection{Non-parallel flux of particles}

\begin{figure}
\begin{minipage}[b]{6.5cm}
\begin{center}
\includegraphics[scale=0.5]{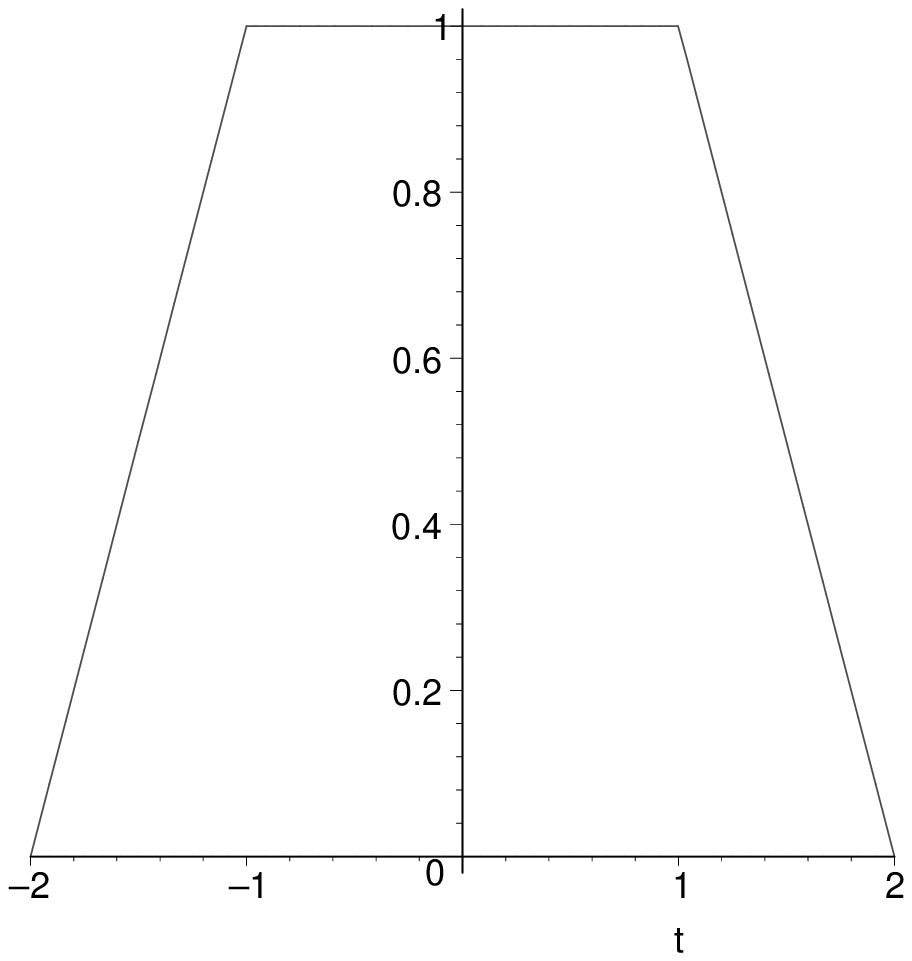}
\end{center}
\caption{$H=1$}
\label{fig:trapezium}
\end{minipage}
\begin{minipage}[b]{6.5cm}
\begin{center}
\includegraphics[scale=0.5]{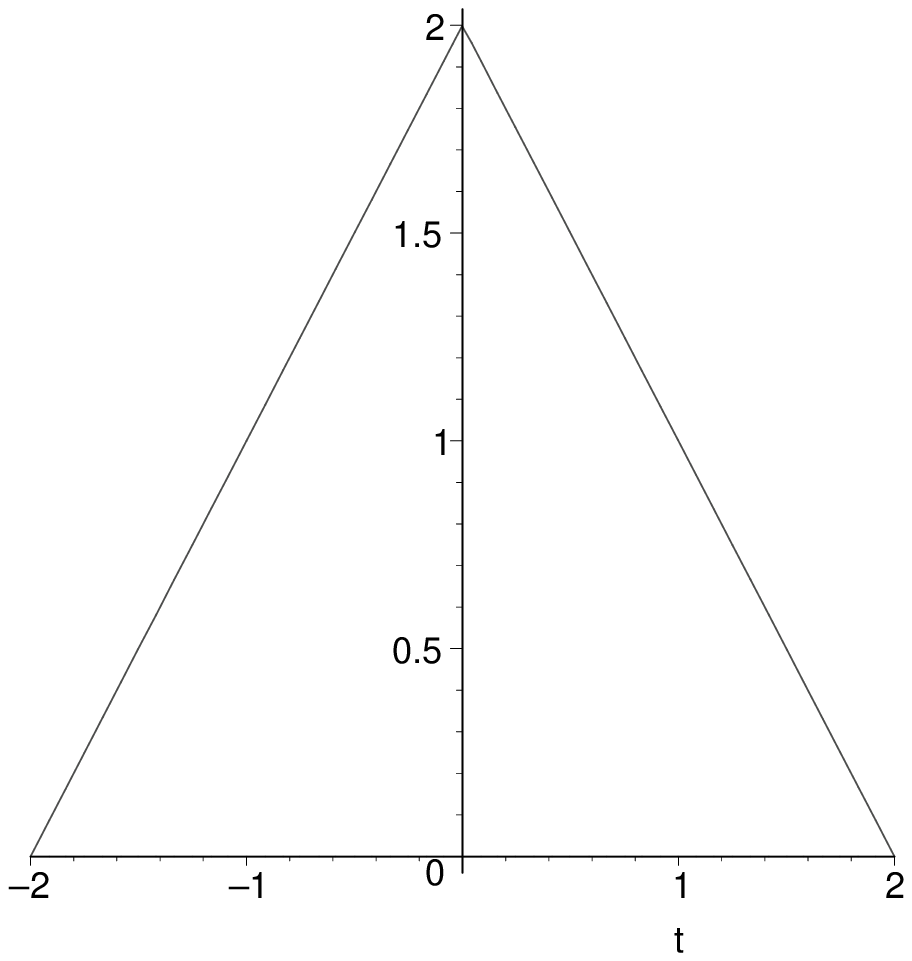}
\end{center}
\caption{$H=2$}
\label{fig:triangle}
\end{minipage}

\bigskip \bigskip

\begin{minipage}[b]{6.5cm}
\begin{center}
\includegraphics[scale=0.5]{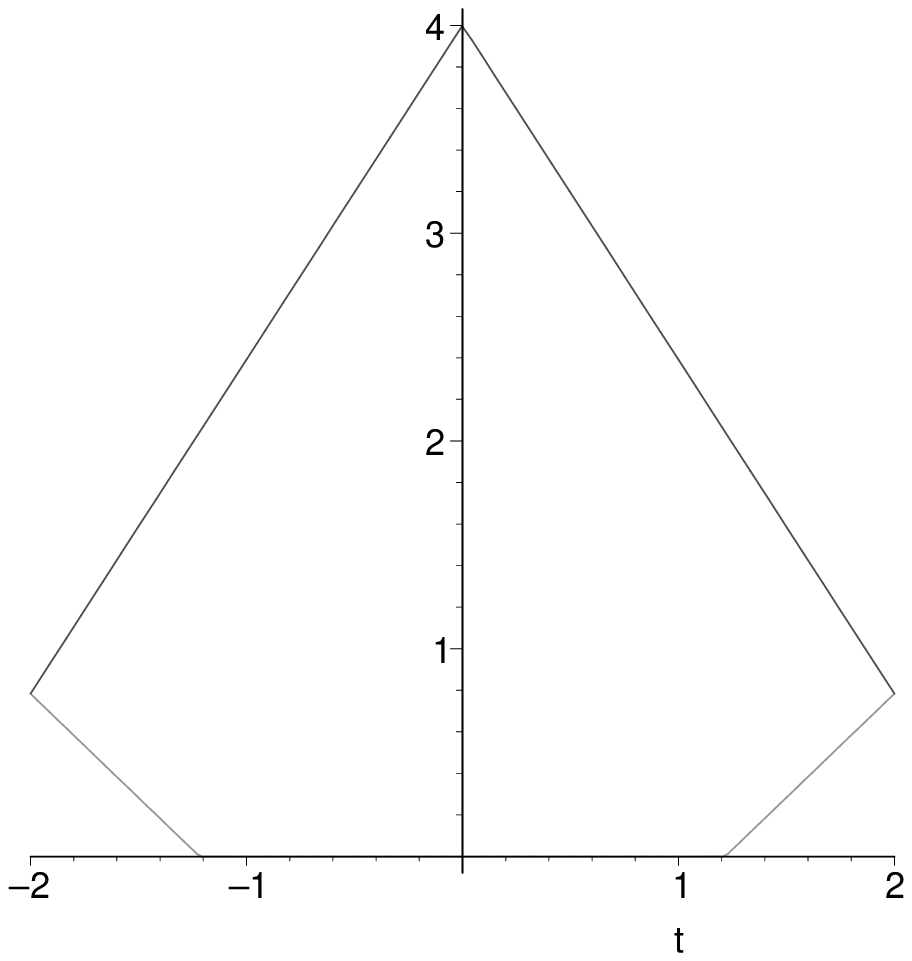}
\end{center}
\caption{$H=4$}
\label{fig:tri-trap}
\end{minipage}
\begin{minipage}[b]{6.5cm}
\begin{center}
\includegraphics[scale=0.5]{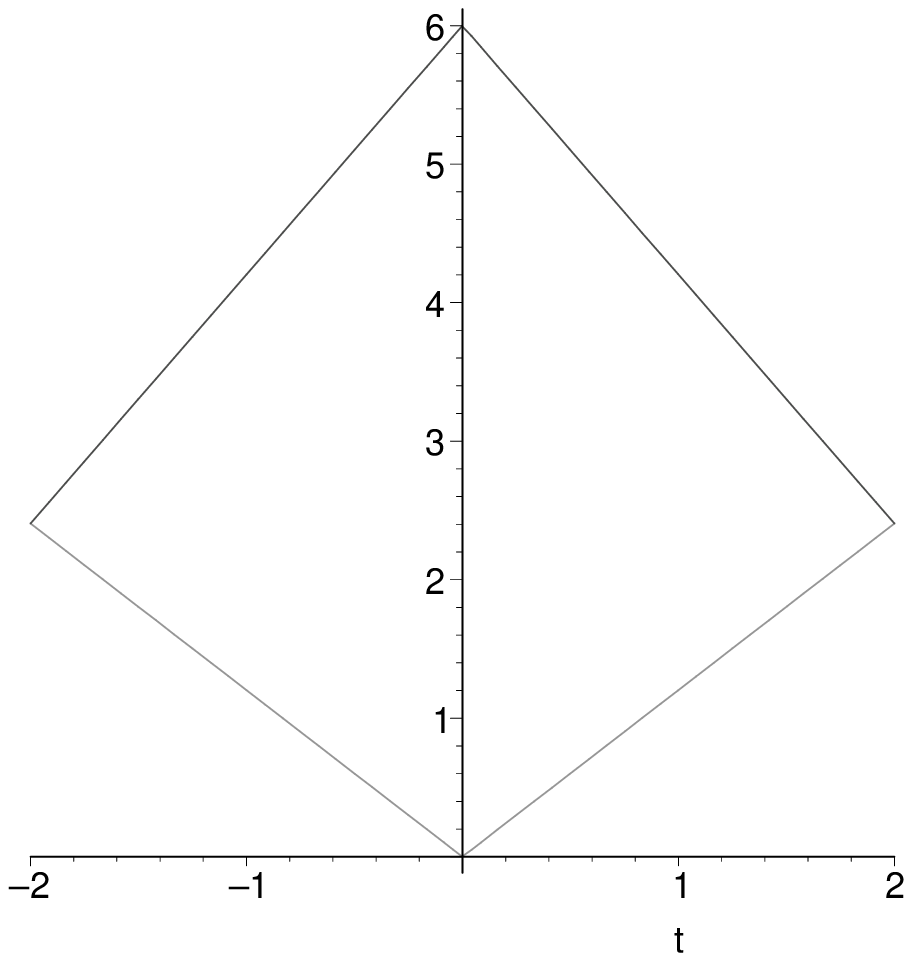}
\end{center}
\caption{$H=6$}
\label{fig:tri-tri}
\end{minipage}

\bigskip \bigskip

Solutions of the two-dimensional Newton-type problem
with $p_+ = \frac{1}{1+u^2} + 0.5$,
$p_- = \frac{0.5}{1+u^2}-0.5$ (non-parallel flux of particles),
$T = 2$, and different values for the height $H$ of the body.
\end{figure}

Let us consider the two-dimensional case ($d = 2$).
As proved in \S \ref{d=2}, there exist four possible cases.
To illustrate this we choose, as an example,
the pressure of the front part of the surface to be
$p_+ = \frac{1}{1+u^2} + 0.5$;
the pressure on the rear part given by
$p_- = \frac{0.5}{1+u^2}-0.5$; the radius $T$
of the maximal cross section of the body to be two
($T = 2$); and then we choose different values for
the height $H$ of the body. Applying the formulas
given in \S \ref{d=2} one obtains that for $H = 1$
the solution is a trapezium (Fig.~\ref{fig:trapezium});
for $H = 2$ a triangle (Fig.~\ref{fig:triangle});
for $H = 4$ the union of a triangle and a trapezium
turned over (Fig.~\ref{fig:tri-trap});
and for $H = 6$ the union of two triangles
with common base (Fig.~\ref{fig:tri-tri}).
We remark that in Newton's problem one has
$p_+ = \frac{1}{1+u^2}$ and $p_- = 0$ (parallel flux),
and only the first two situations occur: solution
to Newton's two-dimensional problem is either
a trapezium or a triangle.

The two-dimensional problem under a non-parallel flux
of particles with density of distribution over velocities
circular gaussian, with biased mean, is studied
in \cite{controlo2004}.


\subsubsection{Newton's classical problem}
\label{Sec:NewPrb}

We now obtain the well-known Newton's solution. For that
we fix $d = 3$, $p_+(u) = 1/(1 + u^2)$, and $p_-(u) = 0$.
Applying the method described in \S \ref{d ge 3},
after some algebra one obtains $\bar{u}_+ =
1/\sqrt{3}$,\, $u^0 = 1$,\, $B_+ = 1/ 2$,\, $\beta = H$, and the
optimal solution $x(t)$ is given in parametric form by
\begin{equation*}
\begin{split}
x &= \frac{\lambda}{2} \left( \frac{3u^4}{4} + u^2 - \ln u - \frac 74
\right) \, , \\
t &= \frac{\lambda}{2} \left( u^3 + 2u + \frac 1u \right) \, ,
\quad 1 \le u \le U \, ,
\end{split}
\end{equation*}
all in agreement with classical formulas.
Expressing the formulas with respect to $U$ and $T$ one obtains:
\begin{gather*}
\lambda = \frac{2 T U}{\left(1 + U^{2} \right)^{2}} \, , \quad
t_0 = \frac{4 T U}{\left(1+{U}^{2}\right)^{2}} \, , \quad
\beta = \frac{T U \left(-7 + 4 U^{2} + 3 U^{4} - 4 \ln(U)\right)}{4\left(1+U^{2}\right)^{2}} \, ,\\
t = \frac{T U\left(1+{u}^{2}\right)^{2}}{u\left(1+U^{2}\right)^{2}} \, , \quad
x = \frac{T U \left(-7 + 4 u^{2}+3 u^{4}-4 \ln(u)\right)}{4\left(1+U^{2}\right)^{2}} \, , \\
{\cal R}_+ = \frac {{T}^{2}\left(17 U^{2}+2+10 U^{4} + 3 U^{6}+
4\,\ln(U) U^{2}\right)}{4 \left(1+U^{2}\right)^{4}} \, .
\end{gather*}
In this case ${\cal R}_- = 0$.


\subsubsection{Newton's problem in higher dimensions}
\label{Sec:NewPrbHD}

Our approach to Newton's problem is valid for an arbitrary
$d \ge 2$. For example, for $d = 4$ (problem in dimension four) one gets:
\begin{gather*}
\lambda = {\frac{2\,{T}^{2} U}{\left (1+{U}^{2}\right )^{2}}} \, , \quad
t_0 = \sqrt {{\frac {4{T}^{2} U}{\left (1+{U}^{2}\right )^{2}}}} \, , \quad
\beta = {\frac {T\left (-5\,U+3\,{U}^{3}+2\sqrt{U}\right )}{5\left(1+{U}^{2}\right)}} \, , \\
t = T \sqrt{\frac{U}{u}} \left( \frac{1+u^2}{1+U^2}\right) \, , \quad
x = {\frac {T\,\sqrt {U}\left (-5\,\sqrt {u}+3\,{u}^{5/2}
+2\right )}{5 \left(1+{U}^{2}\right)}} \, ,\\
{\cal R}_+ = \frac{T^{2}\left(1+3U^{2}\right)}{2 \left(1+U^{2}\right )^{2}} \, .
\end{gather*}


\section*{Acknowledgements}

Research partially supported by the R\&D unit
\emph{Centre for Research in Optimization and Control} (CEOC)
of the University of Aveiro, through the
Portuguese Foundation for Science and Technology (FCT),
cofinanced by the European Community fund FEDER.
The authors also thank an anonymous referee for comments.



\end{document}